\pdfoutput=1 
\documentclass[12pt]{amsart}
\usepackage[T1]{fontenc}
\usepackage[unicode]{hyperref}
\usepackage{color}
\usepackage{amsfonts,amssymb,amsmath,amscd,amsthm}
\usepackage{tikz}
\usetikzlibrary{cd}
\providecommand{\arxiv}[1]{\href{http://arxiv.org/abs/#1}{arXiv:#1}}

\makeatletter
\def\@settitle{\begin{center}%
		\baselineskip14\p@\relax
		\bfseries
		\@title
	\end{center}
}
\makeatother

\theoremstyle{definition}
\newtheorem*{notn}{Notation}
\newtheorem*{defn}{Definition}
\theoremstyle{remark}
\newtheorem*{remark}{Remark}
\newtheorem{case}{Case}
\theoremstyle{plain}
\newtheorem*{theorem}{Theorem}

\newtheorem*{conjecture*}{Conjecture}

\newtheorem*{corol}{Corollary}

\theoremstyle{definition}

\theoremstyle{remark}

\theoremstyle{plain}
\newtheorem{theorem_n}{Theorem}
\newtheorem{propn_n}{Proposition}

\def \k {\mathbf{k}}
\def \F {\mathbf{F}}
\def \A {\mathbb{A}}
\def \P {\mathbb{P}}
\def \fp {\mathfrak{p}}
\def \Supp {\mathrm{Supp}}
\def \Tame {\mathrm{Tame}}
\def \Div {\mathrm{div}}

\def \mult {\mathrm{mult}}

\newcommand{\stein}[2] {\left\lbrace #1,#2 \right\rbrace}
\makeatletter
\newcommand*{\defeq}{\mathrel{\rlap{%
			\raisebox{0.3ex}{$\m@th\cdot$}}%
		\raisebox{-0.3ex}{$\m@th\cdot$}}%
	=}
\makeatother
\newcommand{\relmiddle}[1]{\mathrel{}\middle#1\mathrel{}}
\newcommand{\RN}[1]{
	\textup{\uppercase\expandafter{\romannumeral#1}}
}
\title{Rational endomorphisms of plane preserving a rational volume form}
\author{Georgy Belousov}
\theoremstyle{Definition}
\newtheorem*{jac}{Jacobian conjecture}

\begin{document}
\maketitle
\begin{abstract}
	Let $\varphi$ be a rational map $\P^2 \dashrightarrow\P^2$ that preserves the rational volume form $\frac{\mathrm{d}x}{x}\wedge\frac{\mathrm{d}y}{y}$. Sergey Galkin conjectured that in this case $\varphi$ is necessarily birational. We show that such a map preserves the element $\{x,y\}$ of the second K-group $K_2(\k(x,y))$ up to multiplication by a constant, and restate this condition explicitly in terms of mutual intersections of the divisors of coordinates of $\varphi$ in a way suitable for computations.
	\end{abstract}
\begin{section}{Introduction}
	Diller and Lin in their paper \cite{diller} showed that if a general rational self-map $\varphi : \P^2 \dashrightarrow\P^2$ preserves a rational two-form $\omega$ on the complex plane then, up to a birational change of coordinates, one of the following holds:\\
	1) $-\mathrm{div}(\omega)$ is a smooth cubic curve\\
	2) $\omega=\frac{\mathrm{d}x}{x}\wedge\frac{\mathrm{d}y}{y}$\\
	3) $\omega=\mathrm{d}x\wedge\mathrm{d}y$.\\
	In this work we consider the second case of this classification. Let us call a rational map ${\varphi : \P^2 \dashrightarrow\P^2}$ that preserves the rational volume form $\frac{\mathrm{d}x}{x}\wedge\frac{\mathrm{d}y}{y}$ a \emph{symplectic rational map}.\\ Sergey Galkin conjectured that any symplectic rational map is birational.
	\par We make some preliminary steps towards this conjecture, namely:\\
	1) We show that this conjecture is equivalent to another one, involving $K_2(\k(\P^2))$---see \mbox{Theorem \ref{equiv}}.\\
	2) A description of $K_2$ in terms of algebraic cycles (Proposition \ref{K2/const}) allows us to write down the condition of preserving the rational volume form explicitly in terms of the divisors of coordinates of $\varphi$. The resulting formulas provide a convenient way of checking whether a map is symplectic (see Section \ref{computing}).\\
	3) We provide an example of application of our formulas in the last section of this article.\\
	It may be interesting to compare our conditions with the generators and relations for the group of birational transformations of $\P^2$ preserving $\frac{\mathrm{d}x}{x}\wedge\frac{\mathrm{d}y}{y}$, which were obtained by J\'er\'emy Blanc \cite{blanc}.\\
	\par The author is sincerely grateful to his supervisor Sergey Gorchinskiy for his support and dedication throughout this work, and to Sergey Galkin for suggesting this problem, and for pointing out a mistake in a previous version of the article.
\end{section}
\begin{section}{Statement of the conjecture}
	Let $\k$ be a field of characteristic 0.
	A dominant rational map $\varphi\,$:$\; X \dashrightarrow Y$ of irreducible varieties over $\k$ induces an inverse image on the spaces of rational differential forms over $\k$,\\ denoted by $\varphi^*$:$\;\Omega^{2}_{\k(Y)}\to \Omega^{2}_{\k(X)}$.
	\begin{conjecture*}\label{conj_forms}
		Let $\varphi: \A^2\dashrightarrow\A^2$ be a dominant rational map such that $\varphi^*\left(\dfrac{\mathrm{d}x}{x}\wedge\dfrac{\mathrm{d}y}{y}\right)=\\=\dfrac{\mathrm{d}x}{x}\wedge\dfrac{\mathrm{d}y}{y}\in \Omega^2_{\k(x,y)}$. Then $\varphi$ is birational.
	\end{conjecture*}
	\begin{remark}This conjecture is related to the following famous question.\end{remark}
	\begin{jac}	Let $f: \A^2\rightarrow\A^2$ be a regular map such that $f^*\left(\mathrm{d}x\wedge\mathrm{d}y\right)=\mathrm{d}x\wedge\mathrm{d}y$.\\ Does it follow that $f$ is necessarily biregular?
		\end{jac}
	One possible obstruction to proving the Jacobian conjecture via methods of birational geometry is unconstrained growth of singularities of pullbacks of the form $\mathrm{d}x\wedge\mathrm{d}y$ under birational transformations. These complications do not arise when dealing with the logarithmic form $\dfrac{\mathrm{d}x}{x}\wedge\dfrac{\mathrm{d}y}{y}$.
	There is another version of this conjecture, with Milnor K-groups in place of logarithmic forms. Let us introduce the notation first.
	\begin{defn} (see \cite{kbook} Ch. III \S 7) For a field $F$, let the graded ring $K^{M}_{\bullet}(F)$ be the quotient of the tensor algebra of $F^*$ by the two-sided ideal generated by the elements $\left\{ f\otimes (1-f) \relmiddle{|} 1\neq f\in F^* \right\}$. The \textit{n-th Milnor K-group} $K^{M}_{n}(F)$ is the $n$-th graded component of $K^{M}_{\bullet}(F)$. We will write $\left\lbrace f_1,\dots,f_n \right\rbrace$ for the class of $f_1\otimes\cdots\otimes f_n$ in $K^{M}_{\bullet}(F)$.
	\end{defn}
	\begin{notn}\label{kmodconst}
		Let $F$ be a field extension of $\k$; the group $K^M_2(F)\mathop{/}\left< \stein{F^*}{\k^*} \right>$ will be denoted $K_2(F)/\mathrm{Const}$. For brevity, we will not distinguish between the elements of $K^M_2(F)$ and their images in $K_2(F)/\mathrm{Const}$.
	\end{notn}
	A dominant rational map $\varphi\,$:$\; X \dashrightarrow Y$ of irreducible varieties over $\k$ induces an inverse image $\varphi^*: K^M_\bullet (\k(Y)) \to K^M_\bullet (\k(X))$ on $K$-groups of their function fields. In particular, there is an inverse image on the groups $K_2$ and $K_2/\mathrm{Const}$.
	\begin{conjecture*}\label{conj_k2} 
		Let $\varphi: \A^2 \dashrightarrow \A^2$ be a dominant rational map such that $\varphi^*\stein{x}{y}=\\=\stein{x}{y}\in K_2/\mathrm{Const}$. Then $\varphi$ is birational.
	\end{conjecture*}
	We will show that the two conjectures are in fact equivalent, specifically:
	\begin{theorem_n}\label{equiv}
		A dominant rational map $\varphi: \A^2 \dashrightarrow \A^2$ preserves the form $\frac{dx}{x}\wedge\frac{dy}{y}\in \Omega^2_{\k(x,y)}$ if and only if it preserves the element $\stein{x}{y}\in K_2/\mathrm{Const}$.
	\end{theorem_n}
\end{section}

\begin{section}{$K_2$ in terms of algebraic cycles}
	 For a scheme $X$, let $\mathcal{K}^X_m$ be the Zariski sheaf associated to the presheaf $U \mapsto K_m(U)$.
	\begin{theorem}(Quillen; see \cite{kbook} ch.V Prop 9.8.1)\\Let $X$ be a regular quasi-projective scheme over the field $\k$; then for each $m\geq0$ the complex
		\begin{equation}\label{Gersten}
		0
		\rightarrow \mathcal{K}^X_m
		\xrightarrow{d_0} \bigoplus_{\fp\in X^{(0)}} i_{\fp,*} K_{m}(\kappa(\fp))
		\xrightarrow{d_{1}} \bigoplus_{\fp\in X^{(1)}} i_{\fp,*} K_{m-1}(\kappa(\fp))
		\xrightarrow{d_2}\cdots
		\xrightarrow{d_{m}} \bigoplus_{\fp\in X^{(m)}} i_{\fp,*} \mathbb{Z}
		\rightarrow 0,
		\end{equation} 
		where $\begin{tikzcd}\mathrm{Spec}(\kappa(\fp)) \arrow[r, hook, "\displaystyle{i_\fp}"]& X\end{tikzcd}$ denotes the natural inclusion map, provides a flabby resolution for the $\mathcal{K}_m$-sheaf of $X$.
	\end{theorem}
	Let $X=\A_{\k}^n$, and fix an identification $\F\defeq\k(x_1,...,x_n)\simeq\k(X)$. Applying Gersten's construction we obtain the following flabby resolution for the sheaf $\mathcal{K}^X_2$ on $X$:
	\begin{equation*}
	%\mathcal{K}_2\mid_{\A^2}\rightarrow 
	G \defeq
	\left[
	i_{0,*}K_{2} \left( \F \right)
	\rightarrow	\bigoplus_{\fp\in X^{(1)}} i_{\fp,*}\kappa(\fp)^*
	\rightarrow	\bigoplus_{\fp\in X^{(2)}} i_{\fp,*}\mathbb{Z}
	\rightarrow 0
	\right].
	\end{equation*}
	\\ We can compute the cohomology of $\mathcal{K}^{\A^n}_2$, and therefore of $G$, using $\A^{1}$-invariance of cohomology for the $\mathcal{K}_m$-sheaf (see \cite{rr}).
	\[
	H^i(X,G)=H^i(\A_{\k}^n,\mathcal{K}_2)=
	\begin{cases}
	K_2(\k), & i=0\\
	0, & i\neq 0
	\end{cases}
	\]
	Therefore
	\begin{equation}\label{affinegersten}
	0
	\rightarrow K_2(\k) 
	\xrightarrow{d_0} K_2(\F)
	\xrightarrow{d_1 \iffalse \Tame \fi}	\bigoplus_{\fp\in X^{(1)}} \k{(D_\fp)}^*
	\xrightarrow{d_2 \iffalse \Div \fi}	\bigoplus_{\fp\in X^{(2)}} \mathbb{Z}
	\rightarrow 0
	\end{equation}
	is an exact sequence of abelian groups; here $D_\fp$ stands for the irreducible divisor $V(\fp)\subset X$ defined by the prime ideal $\fp\in X^{(1)}$. The group $\bigoplus_{\fp\in X^{(2)}} \mathbb{Z}$ is none other than the group $Z^2(X)$ of zero-cycles. The differentials of this complex are given by direct sums of residue maps along all irreducible divisors; $d_1$ is called the \emph{tame symbol}, and can be computed (see \cite{kbook} Ch.$\:$III Lemma 6.3) as
	\begin{equation}\label{tame}\Tame\stein{x}{y}=\bigoplus_{\fp \in X^{(1)}} (-1)^{\nu_\fp(x) \nu_\fp(y)}\dfrac{y^{\nu_\fp(x)}}{x^{\nu_\fp(y)}}(\mathrm{mod}\, \fp).\end{equation}
	Each component of $d_2$ maps a function to the corresponding principal divisor:
	\[\k(D)^* \xrightarrow{\Div} Z^2(D) \subset Z^2(X). \]
	\\We will be interested in describing the group $K_2(\F)/\mathrm{Const}$; notice that
	\begin{equation}\label{ker(Tame)} d_0\left(K_2(\k)\right) \subset \stein{\F^*}{\k^*} \subset K_2(\F).\end{equation}
	Thus the exact sequence (\ref{affinegersten}) yields
	\[
	K_2(\F)\mathop{/} \left<\stein{\F^*}{\k^*}\right> \simeq 
	\left.\ker\left(
	\bigoplus_{D \subset \A^n} \k(D)^*
	\xrightarrow{\Div}	Z^2(\A^n)
	\right)
	\right/ \Tame\left<\stein{\F^*}{\k^*}\right>.
	\]
	The image of the group $\mathrm{Const}$ is just the set of sums of constant functions, as shown in the following proposition.
	\begin{propn_n}\label{Tame(const)}
		\[\Tame\left<\stein{\F^*}{\k^*}\right> 
		= \bigoplus_{\fp\in X^{(1)}}\k^*
		\hookrightarrow \bigoplus_{\fp\in X^{(1)}}\k(D_\fp)^* \text{ (included as sums of constant functions.)}\]
	\end{propn_n}
	\begin{proof}
		We have
		\[ \Tame \stein{f}{c} = \bigoplus_{\fp} \left(\dfrac{c^{\nu_{\fp}(f)}}{f^0}\right)_{\mathrm{mod}\:\fp} =\: \bigoplus_{\fp} c^{\nu_{\fp}(f)}\]
		Since the divisor class group is trivial for $X=\A^n$, for any $\fp\in X^{(1)}$ we can choose an equation $F_\fp$ defining the corresponding divisor $D_\fp$. Then for any finite collection of $\fp_i$ and $c_i$ we can take $\psi=\sum \stein{F_{\fp_i}}{c_i}$ in $K_2(\F)$ so that $\Tame(\psi)=\bigoplus_{\fp_i} c_i$; this shows that any finite collection of constant functions is contained in $\Tame\left<\stein{\F^*}{\k^*}\right>$.
	\end{proof}
	This gives our final presentation for $K_2(\F)/\mathrm{Const}$, concluding this section.
	\begin{propn_n}\label{K2/const}$\displaystyle{
		K_2(\F)\mathop{/} \mathrm{Const} \simeq 
		\ker\left(
		\bigoplus_{D\subset \A^n} \k(D)^*\mathop{/}\k^*
		\xrightarrow{\Div}	Z^2(\A^n)
		\right)}.$
	\end{propn_n}
\end{section}

\begin{section}{Logarithmic differential and proof of theorem 1}
	We will need the notion of rational logarithmic form (see \cite{polar}, \S 2 for a thorough exposition).
	\begin{defn} 
		Let $F$ be an extension of $\k$ of finite transcendence degree. A rational differential form $\omega \in \Omega^p_{F/\k}$ is called \emph{logarithmic} if there exists a smooth projective variety $V/\k$ with function field $F$ and a simple normal crossing divisor $D\subset V$, such that $\omega$ is locally logarithmic along $D$ and regular elsewhere on $V$. We will write $\Omega^n_{F,\, log}$ for the subgroup of $\Omega^n_F$ generated by logarithmic forms.
	\end{defn}
	There is a map $K_n^M(F) \xrightarrow{\mathrm{dlog}} \Omega^n_F$ taking $\left\{f_1, \dots, f_n\right\}$ to the form $\frac{\mathrm{d}f_1}{f_1}\wedge\dots\wedge\frac{\mathrm{d}f_n}{f_n}$. Such a form is indeed logarithmic.
	\begin{theorem_n} The map$\begin{tikzcd}K_2(\F)/\mathrm{Const} \arrow[r, "{\mathrm{dlog}}"]&\Omega^2_{\F,\,log}\end{tikzcd}$ is injective.
	\end{theorem_n}
	\begin{proof}
				The following diagram is commutative:
				\[\begin{tikzcd}[%
				nodes={asymmetrical rectangle}
				,/tikz/column 1/.append style={anchor=base east}
				,/tikz/column 2/.append style={anchor=base west}
				]
				K_2(\mathbf{F})/\mathrm{Const} \arrow[r, hook, "\Tame"] \arrow[d, start anchor = south east, end anchor = north east, shorten <=0.8em, shift right = 1.8 em, pos=.63,"{\mathrm{dlog}}"]&[1em]
				\displaystyle \bigoplus_{D \subset X} \k(D)^*/\k^* \arrow[d, hook, start anchor = south west, end anchor = north west, shift left = 1.3 em, "{\mathrm{dlog}}"]
				\\[0.5em]
				\Omega^2_{\mathbf{F},\, log} \arrow[r, "{\mathrm{Res}}"] &
				\displaystyle \bigoplus_{D \subset X} \Omega^1_{\k(D),\, log}
				\end{tikzcd}\]
				(One can check this by applying $\mathrm{dlog}$ to the formula (\ref{tame}) for $\Tame$.)\\
				The top horizontal map $\Tame$ is injective by Proposition \ref{K2/const}.\\
				The right vertical map $\mathrm{dlog}$ is injective, as $\mathrm{dlog}$ of a function vanishes if and only if the function is constant. Therefore the left vertical map is also injective.
	\end{proof}
		Recall that our goal was to prove that a rational map $\varphi$$\,:\A^2 \dashrightarrow \A^2$ preserves $\stein{x}{y}\in K_2/\mathrm{Const}$ whenever it preserves the logarithmic volume form $\frac{\mathrm{d}x}{x}\wedge\frac{\mathrm{d}y}{y}\in\Omega^{2}_{\k(x,y),\:log}$ (Theorem \ref{equiv}).
	\begin{proof}[Proof of Theorem 1]
		If $\varphi$ leaves $\stein{x}{y}\in K_2(\k(x,y))/\mathrm{Const}$ in place, then it clearly preserves $\frac{\mathrm{d}x}{x}\wedge\frac{\mathrm{d}y}{y}=\mathrm{dlog}\stein{x}{y}$, as $\varphi^*$ commutes with $\mathrm{dlog}$.\\
		Suppose $\varphi$ preserves $\frac{\mathrm{d}x}{x}\wedge\frac{\mathrm{d}y}{y}$. Then it also preserves $\mathrm{dlog}^{-1}\left(\frac{\mathrm{d}x}{x}\wedge\frac{\mathrm{d}y}{y}\right)$ as a set.
		But since we have observed that $\mathrm{dlog}$ is an inclusion, $\stein{x}{y}$ is the only element of $\mathrm{dlog}^{-1}\left(\frac{\mathrm{d}x}{x}\wedge\frac{\mathrm{d}y}{y}\right)$; this finishes the proof.
	\end{proof}
\end{section}
\begin{section}{Computing tame symbol}\label{computing}
	Let us fix the projectivization $\overline{\A^2}=\P^2$. Given a divisor $D$ on $\A^2$, denote its closure in $\P^2$ by $\overline{D}$ for a moment. The divisorial map $\k(D)^*/\k^*\rightarrow\mathrm{PDiv}(\overline{D})$ is an isomorphism, which implies the following.
	\begin{equation}\label{description}
		K_2(\k(x,y))/\mathrm{Const} \simeq \ker \left( \bigoplus_{D \subset \A^2}\mathrm{PDiv}\left(\overline{D}\right)\xrightarrow{\Div} Z^2(\A^2)\right)
	\end{equation}
	 The last proposition can be rephrased to obtain the following combinatorial description of the group on the left: an element of $K_2(\k(x,y))/\mathrm{Const}$ amounts to a finite collection of divisors $D_i$ on $\A^2$, and a choice of a principal divisor $\alpha_i$ on the projectivization $\overline{D_i}$ of each of the $D_i$, such that $\sum_i \mult_x(\alpha_i) = 0$ for all closed points $x\in\A^2$. Notice that we do not impose any conditions on $\alpha_i$ at the boundary $\P^2\setminus\A^2$.\\
	This description of elements of $K_2(\k(x,y))/\mathrm{Const}$ is useful for explicit computations with rational maps. Since it involves the closures of divisors, we will now switch to considering principal divisors on $\P^2$ in place of affine divisors for brevity.
	\begin{notn}
		1) For a pair of irreducible divisors $D,E$ on $X$, \\define $D\sqcap E\in Z^2(X)$, ${D \sqcap E}\defeq\begin{cases}
		D \cap E, & \text{if }D \neq E\\
		0, & \text{if } D=E
		\end{cases}$\\
		This operation extends linearly over the whole $\mathrm{Div}(X)$.\\
		2) There exist functions $f_D, f_E \in \k(\A^2)$ such that $D=(f_D), E=(f_E)$. Since these functions are determined uniquely up to a constant multiplier, $\Tame\stein{D}{E}\defeq \Tame\stein{f_D}{f_E} \in K_2(\k(x,y))/\mathrm{Const}$ is well-defined.
	\end{notn}
	Given a pair of principal divisors $D,E$ on $\P^2$, let us express the components of $\Tame\stein{D}{E} \in K_2(\k(x,y))/\mathrm{Const}$ in terms of the original divisors. If $\displaystyle{\Tame\stein{D}{E} = \bigoplus_{l_\infty\neq C \subset \P^2} \alpha_C}$, then $\alpha_C$ can be computed as follows:
	\begin{flalign*}
	&\text{If }C \notin \mathrm{Supp}(D) \cup \mathrm{Supp}(E)\text{, then } \alpha_C = 0.\\
	&\text{Case \RN{1}:}\text{ if }C \in \mathrm{Supp}(D) \setminus \mathrm{Supp}(E), \\
	&\;\;\;\text{then }\alpha_C=\nu_C(D)\cdot\left( C \sqcap E \right).\\
	&\text{Case \RN{2}:}\text{ if }C \in \mathrm{Supp}(E) \setminus \mathrm{Supp}(D), \\
	&\;\;\;\text{then }\alpha_C=-\nu_C(E)\cdot\left( C \sqcap D \right).\\
	&\text{Case \RN{3}:}\text{ if }C \in \mathrm{Supp}(D) \cap \mathrm{Supp}(E), \\
	&\;\;\;\text{then }\alpha_C=\nu_C(D)\cdot\left( C \sqcap E \right) - \nu_C(E)\cdot\left( C \sqcap D \right).&
	\end{flalign*}
	Now suppose $\Tame\stein{D}{E}=\Tame\stein{x}{y}$. This is equivalent to $\alpha_C$ being equal to the $C$-component of $\Tame\stein{x}{y}$ for all divisors $l_\infty \neq C\subset \P^2$. 
	\begin{equation}\label{tamexy}
		\Tame\stein{x}{y}=\underbrace{(0-\infty_h)}_{l_h}\oplus\underbrace{(-0+\infty_v)}_{l_v},\end{equation}$ \textit{ where } l_h=\{y=0\}\text{, }l_v=\{x=0\}\text{, }0=[0:0:1]\text{, }\infty_h=[1:0:0]\text{, }\infty_v=[0:1:0].
	$\\
	All the components of $\Tame\stein{x}{y}$ except $l_h$ and $l_v$ are trivial, i.e. $\forall C \neq l_h, l_v:\; \alpha_C = 0$, which has the following implications for $D$ and $E\,$:\\
	1) $C \in \mathrm{Supp}(D)$, $C \notin \mathrm{Supp}(E) \cup \lbrace l_h, l_v \rbrace$, \\
	then $0=\alpha_C=\nu_C(E)\cdot(C \sqcap D)$ $\Leftrightarrow$ $\nu_C(D)=1$, $Div^0(C) \ni C\sqcap E = 0$.\\
	2) $C \in \mathrm{Supp}(E)$, $C \notin \mathrm{Supp}(D) \cup \lbrace l_h, l_v \rbrace$, \\
	then $0=\alpha_C=\nu_C(D)\cdot(C \sqcap E)$ $\Leftrightarrow$ $\nu_C(E)=1$, $Div^0(C) \ni C\sqcap D = 0$.\\
	3) $C \in \mathrm{Supp}(E) \cup \mathrm{Supp}(D)$, $C \notin \lbrace l_h, l_v \rbrace$, \\
	then $0=\alpha_C=\nu_C(D) \cdot (C \sqcap E )) - \nu_C(E) \cdot (C \sqcap D))$\\
	For $l_h$ and $l_v$, we have similar conditions depending on which of the divisors $D$ and $E$ contain these lines as their components. From symmetry considerations there are four distinct cases:\\
	
	Case $\left(\RN{1},\RN{1}\right)$:\\
	$\nu_h(D)=\pm1, \nu_v(D)=1$\\ $\nu_h(E)=\nu_v(E)=0$\\
	$l_h\sqcap E = \pm(0 - \infty_h)$\\
	$l_v\sqcap E = -0 + \infty_v$\\
		
	Case $\left(\RN{1},\RN{2}\right)$:\\
	$\nu_h(D)=1, \nu_v(D)=0$\\ $\nu_h(E)=0, \nu_v(E)=1$\\
	$l_h\sqcap E = 0 - \infty_h$\\
	$l_v\sqcap D = -0 + \infty_v$\\
	
	Case $\left(\RN{1},\RN{3}\right)$:\\
	$\nu_h(D)=\pm1, \nu_v(D)\geq 1$\\ $\nu_h(E)=0, \nu_v(E)\geq 1$\\
	$0 - \infty_h = \alpha_{l_h} =\pm(l_h\sqcap E)$\\
	$-0+\infty_v=\alpha_{l_v}=\nu_v(D)\cdot\left( {l_v} \sqcap E \right) - \nu_v(E)\cdot\left( {l_v} \sqcap D \right)$\\
	
	Case $\left(\RN{3},\RN{3}\right)$:\\
	$\left|\nu_h(D)\right|\geq 1, \nu_v(D)\geq 1$\\ $\nu_h(E)\geq 1, \nu_v(E)\geq 1$\\
	$0-\infty_h=\alpha_{l_h}=\nu_h(D)\cdot ( l_h \sqcap E ) - \nu_h(E)\cdot ( l_h \sqcap D)$\\
	$-0+\infty_v=\alpha_{l_v}=\nu_v(D)\cdot ( l_v \sqcap E ) - \nu_v(E)\cdot ( l_v \sqcap D)$.
\end{section}
\begin{section}{An example of application}
	Let $\varphi:\P^2 \dashrightarrow \P^2$ be a rational map given by its affine coordinates $\varphi=(f,g)$. Denote by $D,E\in \mathrm{PDiv}(\P^2)$ the corresponding principal divisors, and let $D=D^+ - D^-$, where $D^+$ and $D^-$ are effective divisors without common components. Define $E^+, E^-$ in the same way.\\
	In what follows $\mathrm{Bs}(f), \mathrm{Bs}(g)$ will refer to the indeterminacy loci of $f$, $g$ considered as rational maps $\P^2\dashrightarrow\P^1$. Note that $\mathrm{Bs}(f)=|D^+| \cap |D^-|$, $\mathrm{Bs}(g)=|E^+| \cap |E^-|$. For $C$ an irreducible divisor on $\P^2$, $\alpha_C$ stands for the $C$-component of $\Tame\stein{D}{E}$, as in section \ref{computing}.
	\begin{remark}\label{baseloc_bir}
		Let $\varphi:\P^2 \dashrightarrow \P^2$ be a rational map given by its affine coordinates $\varphi=(f,g)$. If $\mathrm{Bs}(f)\cup\mathrm{Bs}(g) \cup ( D^- \cap E^- ) =\varnothing$, then $\varphi$ is regular.
	\end{remark}
	\begin{proof}
		Let $\varphi = [F:G:H]$, where $F, G, H$ are homogeneous polynomials of the same degree; we may assume they have no common multiples.
		For brevity, let us denote the corresponding \emph{effective} divisors on $\P^2$, as well as their supports, by the same letters $F,G,H$.
		Let 
		\[
			F = D^+ + \Delta_F,
		\]
		\[
			G = E^+ + \Delta_G,
		\]
		\[
			H = D^- + \Delta_F = E^- + \Delta_G,
		\]
		where $D^+$, $D^-$, $E^+$, $E^-$, $\Delta_F$, $\Delta_G$ are effective divisors, $D^+$ has no common components with $D^-$, $E^+$ has no common components with $E^-$, and $\Delta_F$ has no common components with $\Delta_G$. Here, as before, $(f)=D^+ - D^-$, $(g)=E^+ - E^-$. Then we have:
		\begin{align*}
			F &\subset (F - \Delta_F) \cup (H - \Delta_G) = D^+ \cup E^-, \\
			G &\subset (G - \Delta_G) \cup (H - \Delta_F) = E^+ \cup D^-, \\
			H &\subset (H - \Delta_F) \cup (H - \Delta_G) = D^- \cup E^-. 
		\end{align*}
		Thus
		\[	
			Bs(\varphi) = F \cap G \cap H \subset (D^+ \cup E^-) \cap (E^+ \cup D^-) \cap (D^- \cup E^-) = (D^+ \cap D^-) \cup (E^+ \cap E^-) \cup (D^- \cap E^-).
		\]
	\end{proof}
	\begin{propn_n}\label{baseloc}
		Let $(f,g)=\varphi:\P^2 \dashrightarrow \P^2$ be a rational map that preserves $\Tame\stein{x}{y}$. Without loss of generality, suppose that $l_h\in\Supp(D)$, $l_v\in\Supp(E)$. Assume further that $D$ and $E$ have no common components except $l_\infty$ (so that we are in the situation of Case $\left(\RN{1},\RN{2}\right)$). If no three out of the four divisors $D^+, D^-, E^+, E^-$ have nonempty intersection, then $D=l_h-l_\infty$, $E=l_v-l_\infty$.
	\end{propn_n}
	\begin{proof}
		Let $C_1$ be a simple component of $D$, $C_2$ a simple component of $E$. By B\'ezout's theorem, their intersection $\left|C_1\sqcap C_2\right|$ is nonempty. Pick a point $p \in \left|C_1\sqcap C_2\right|$.
		Since no three out of the four divisors $D^+, D^-, E^+, E^-$ have nonempty intersection, the point $p$ is contained in the support of exactly two of these four divisors. 
		\begin{case} 
			Suppose $C_1\neq l_\infty$ and $C_2 \neq l_\infty$. Since $D$ and $E$ have no common components except possibly $l_\infty$, both $\alpha_{C_1}=\nu_{C_1}(D) \cdot (C_1 \sqcap E)$ and  $\alpha_{C_2}= - \nu_{C_2}(E) \cdot (C_2 \sqcap D)$ must be nontrivial. Since the only nontrivial components of $\Tame\stein{D}{E}=\Tame\stein{x}{y}$ are $\alpha_{l_h}$ and $\alpha_{l_v}$, we must have either $C_1 = l_v$ and $C_2=l_h$ or $C_1 = l_h$ and $C_2 = l_v$, and in both cases $p=0 \in \P^2 \setminus l_\infty$.
		\end{case}
		\begin{case}
			Suppose $C_1$ coincides with $l_\infty$. Then $C_2 \neq l_\infty$, and $\alpha_{C_2}= - \nu_{C_2}(E) \cdot (C_2 \sqcap D) \neq 0$. Therefore we must have either $C_2 = l_h$ or $C_2 = l_v$.
			Analogously, if $C_2$ coincides with $l_\infty$, then $C_1$ must coincide with either $l_h$ or $l_v$.
		\end{case}
		Therefore the only possibilities for simple components of $D$ and $E$ are $l_v$, $l_h$ and $l_\infty$. Since we are in the situation of Case $(\RN{1}, \RN{2})$ above, it follows that $D=l_h - l_\infty$, $E = l_v - l_\infty$.
	\end{proof}
	\begin{corol}\label{exapp}
		Let $(f,g)=\varphi:\P^2 \dashrightarrow \P^2$ be a rational map that preserves the form $\dfrac{\mathrm{d}x}{x}\wedge\dfrac{\mathrm{d}y}{y}$, and assume that the divisors $(f)$ and $(g)$ have no common components except possibly $l_\infty$.\\
		If no three out of the four divisors $D^+, D^-, E^+, E^-$ have nonempty intersection, then $\varphi$ is biregular.
	\end{corol}
		Notice that the condition $\mathrm{Bs}(f)\cup\mathrm{Bs}(g)\cup ( D^- \cap E^- )=\varnothing$ implies that no three out of the four divisors $D^+, D^-, E^+, E^-$ have nonempty intersection.
\end{section}

\end{document}